\documentclass[12pt]{article}
\usepackage{amsmath,amssymb,amsthm, epsfig}
\usepackage{amsfonts}
\usepackage{amsmath}
\usepackage{amssymb}
\usepackage{color}
\usepackage[mathscr]{eucal}

\title{\bf A limiting free boundary problem ruled by Aronsson's equation}

\author{Julio D. Rossi\footnote{Departamento de
Matem\'{a}tica,  FCEyN UBA (1428), Buenos Aires, Argentina.
E-mail: \texttt{ jrossi@dm.uba.ar}} \quad $\&$ \quad Eduardo V.
Teixeira \footnote{Universidade Federal do Cear{\'a}. Departamento
de Matem{\'a}tica. Campus do Pici - Bloco 914. Fortaleza, CE -
Brazil 60.455-760. E-mail: \texttt{eteixeira@ufc.br}}}

\newlength{\hchng}
\newlength{\vchng}
\def \dist {\mathrm{dist}}

\def \suchthat {\ \big | \ }
\def \Lip {\mathrm{Lip}}

\def \J {\mathfrak{J}}

\def \Leb {\mathscr{L}^n}

\newtheorem{theorem}{Theorem}
\newtheorem{lemma}{Lemma}

\theoremstyle{definition}

\theoremstyle{remark}
\newtheorem{remark}{Remark}
\numberwithin{equation}{section}

\newtheorem{lema}{Lemma}[section]
\newtheorem{defi}{Definition}[section]

\newcommand{\intav}[1]{\mathchoice {\mathop{\vrule width 6pt height 3 pt depth  -2.5pt
\kern -8pt \intop}\nolimits_{\kern -6pt#1}} {\mathop{\vrule width
5pt height 3  pt depth -2.6pt \kern -6pt \intop}\nolimits_{#1}}
{\mathop{\vrule width 5pt height 3 pt depth -2.6pt \kern -6pt
\intop}\nolimits_{#1}} {\mathop{\vrule width 5pt height 3 pt depth
-2.6pt \kern -6pt \intop}\nolimits_{#1}}}

\begin{document}
\maketitle

\begin{abstract} We study the behavior of $p$-Dirichlet optimal
design problem with volume constraint for $p$ large. As the limit
as $p$ goes to infinity, we find a limiting free boundary problem
governed by the infinity-Laplacian operator. We establish a
necessary and sufficient condition for uniqueness of the limiting
problem and, under such a condition, we determine precisely the
optimal configuration for the limiting problem. Finally, we
establish convergence results for the free boundaries.

\medskip
\noindent \textbf{Keywords:} Optimal design, Free boundary
problems, Infinite Laplacian.

\medskip
\noindent \textbf{2000 MSC} \textit{Primary}: 35R35, 35J70, 62K05,
49L25.
\end{abstract}


\section{Introduction}

Let $\Omega$ be a smooth bounded domain in the Euclidean space
$\mathbb{R}^n$ and $\alpha$ a fixed positive number less than the
Lebesgue measure of $\Omega$. An optimal design problem with
volume constraint can be generally written as:
\begin{equation}\label{General Optimization Prob. INTROC}
    \textrm{Min } \left \{ \J(\mathcal{O}) \suchthat \mathcal{O} \subset \Omega
    \quad \textrm{ and
    } \quad
    \Leb (\mathcal{O}) \le \alpha \right \}.
\end{equation}
For most of applications, $\J(\mathcal{O})$ has an integral
representation involving functions which are linked to the
competing configuration $\mathcal{O}$ by a prescribed PDE.
\par
The modern history of this line of research probably starts at the
pioneering work of Aguilera, Alt and Caffarelli, \cite{AAC}. In
that paper, the authors address the question of minimizing the
Dirichlet integral when prescribed the volume of the zero set.
Lederman in \cite{Led} establishes similar results for
non-homogeneous minimization problem: $\int |Du|^2 - gu$. Alt,
Caffarelli and Spruck, \cite{ACS}, considered the minimization
problem (\ref{General Optimization Prob. INTROC}) for
$\J(\mathcal{O}) = \int_\Omega \Delta u dX$, where $u$ is the
harmonic function in $\mathcal{O}$, taking a prescribed boundary
data $\varphi$ on $\partial \Omega$ and zero on $\partial
\mathcal{O}$. This is a model for an optimal shape problem in heat
conduction theory with non-constant temperature distribution.
Nonlinear optimal design problems with non-constant temperature
distribution was treated in \cite{Teix}.  The common feature of
the aforementioned works is that all of them are governed by the
Laplacian operator. Their fine analysis rely on the revolutionary
work of Alt and Caffarelli, \cite{AC}.
\par
Just recently the study of optimal design problems ruled by
degenerate quasilinear operators was successfully developed. This
theory is the starting point for the main goal of this present
work which we describe now. Let us consider the problem of
minimizing the $p$-Dirichlet integral with a given positive
boundary data $f$ and prescribed the maximum volume of the
support. More precisely, let us consider the following free
boundary optimization problem:
\begin{equation}\label{problem.p}\tag{$\mathfrak{P}_p$}
    \min \left \{\int_{\Omega} |\nabla u (X)|^p dX \suchthat
    u \in W^{1,p} (\Omega), \ u = f \mbox{ on } \partial \Omega, \
    \Leb \left ( \{u>0\} \right ) \leq \alpha \right\}.
\end{equation}
Existence of a minimizer as well as smoothness properties of its
free boundary have been established in \cite{FMW06} and \cite{OT}.
Further generalizations are addressed in \cite{Teix-Preprint,
Teix-Prep}. In the present paper we are interested in the
asymptotic behavior, as $p$ goes to infinity, of optimal shapes to
problem (\ref{problem.p}). Analytical and geometric features of a
limiting free boundary revels asymptotic information upon optimal
design problem \eqref{problem.p}. Driven by classical
considerations, we are led to consider the following ``limiting
problem":
\begin{equation}\label{limit.problem}\tag{$\mathfrak{P}_\infty$}
\min \left \{ \Lip(u) \suchthat u \in W^{1,\infty}(\Omega), \ u =
f, \ \mbox{on} \ \partial \Omega, \ \Leb \left ( \{u>0\} \right )
\leq \alpha \right\}.
\end{equation}
where $\Lip(u)$ is the Lipschitz constant of $u$:
\begin{equation}\label{Lips.const}
    \Lip(u) = \sup_{x,y} \frac{|u(x) - u(y)|}{|x-y|}.
\end{equation}
Our first concern is to prove that any sequence of minimizers
$u_p$ to problem \eqref{problem.p} converges (up to a subsequence)
to a solution, $u_\infty$, of the limiting problem
\eqref{limit.problem}. In addition, we are interested in finding
the EDP $u_\infty$ satisfies in its set of positivity.  In this
direction and enforcing the fact that $u_\infty$ is an extremal
for the Lipschitz minimization problem, we show that $u_\infty$ is
indeed an absolutely minimizer for the Lipschitz constant within
its set of positively, $\Omega_\infty := \{u_\infty > 0 \}$. That
is, it minimizes the Lipschitz constant in every subdomain of
$\Omega_\infty$ when testing against functions with the same
boundary data, see \cite{ACJ}. Hence it is an $\infty-$harmonic
function in its positivity set. These information are the contents
of the first Theorem in this paper which we state now.
\begin{theorem}\label{theorem 1 - convergence}
    Let $u_p$ be a minimizer of \eqref{problem.p}, then, up to a
    subsequence,
    $$
        u_p \to u_\infty, \qquad  \mbox{ as } p \to \infty,
    $$
    uniformly in $\overline{\Omega}$ and weakly in every
    $W^{1,q} (\Omega)$ for $1<q<\infty$, where $v_\infty$ is a
    minimizer of \eqref{limit.problem}. The limiting function
    $u_\infty$ satisfies the PDE,
    $\Delta_\infty u_\infty = 0$,  in $\{u_\infty >0 \}$
    in the viscosity sense.
    Here $\Delta_\infty u := Du D^2 u (Du)^\text{t}$ is the infinity Laplacian.
\end{theorem}
\par
It is known that under the assumptions $\Omega$ convex and $f
\equiv \text{const.}$, one can prove uniqueness for problem
\eqref{problem.p}, \cite{T3} (see also \cite{AM, HS, HS1, HS2, LV}
for related Bernoulli-type problems). However, uniqueness is not
expected in general for problem \eqref{problem.p}. Surprisingly
enough, under a mild compatibility condition upon $\Lip(f)$,
$\Omega$ and $\alpha$, that does not involve convexity assumption
on $\Omega$, we prove uniqueness for the limiting problem
\eqref{limit.problem}. In particular, any sequence of solutions to
problem \eqref{problem.p} converges to a same optimal limiting
configuration. Such a result such be read as an ``asymptotic
uniqueness phenomenon" for problem \eqref{problem.p}. In addition,
we have precisely found the optimal shape for the limiting problem
\eqref{limit.problem}, that is, it revels where and how optimal
configurations $\Omega_p := \{u_p > 0 \}$ stabilize (see also
Remark \ref{rmk applied}).
\par
More precisely, for our next Theorem we shall work under the
following geometric compatibility condition:
\begin{equation}\label{Comp. Condition - thm 2}\tag{H}
    \Leb \left( \bigcup_{y \in\partial \Omega} B_{\frac{f(y)}{\Lip (f)}}
    (y) \cap \Omega \right) \geq \alpha.
\end{equation}
It is understood that if $f$ is constant, then \eqref{Comp.
Condition - thm 2} is automatically satisfied.
\begin{theorem}\label{theorem2 - uniqueness} Assume \eqref{Comp. Condition - thm 2} and let
$\lambda^{\star}$ be the unique positive real number such that the
domain
$$
    \Omega^\star := \bigcup\limits_{x\in \partial \Omega}
    B_{\frac{f(x)}{\lambda^{\star}}}(x) \cap \Omega
$$
has Lebesgue measure precisely $\alpha$. Then the function
$u_\infty$, defined as
$$
    \left \{
        \begin{array}{rll}
            \Delta_{\infty} u_\infty &=& 0 \qquad \text{ in } \Omega^\star, \\
            u_\infty & = & f \qquad \text{ on } \partial \Omega,  \\
            u_\infty & = & 0 \qquad \text{ on } \partial \Omega^\star
            \cap \Omega,
        \end{array}
    \right.
$$
is the unique minimizer for problem \eqref{limit.problem}. Hence,
if $u_p$ is a minimizer of \eqref{problem.p}, then the whole
sequence $u_p$ converges, $
    u_p \to u_\infty$ ,uniformly in $\overline{\Omega}$ and weakly in every $W^{1,q}
(\Omega)$ for $1<q<\infty$. In addition, $u_\infty$ is given by
the formula,
$$
u_\infty =  \max\limits_{y \in \partial \Omega} \left (f(y)
        - \lambda^{\star} |x-y| \right )_{+}.
$$
\end{theorem}

\begin{remark} As mentioned above,
Theorem \ref{theorem2 - uniqueness} applies in particular to an
important physical situation, namely heat conduction problems with
evenly heated domains, i.e., $f\equiv T$ (constant).
\end{remark}

\begin{remark} \label{rmk applied} From the applied point of view, Theorem \ref{theorem2
- uniqueness} provides a rigorous mathematical proof for the
empirical, and widely employed, intuition that says that the
configuration  $\Omega^\star$ should be approximately an optimal
way of insulating a given body $\Omega$ with temperature
distribution $f$.
\end{remark}

\begin{remark} The fact that the equation that rules the limit configuration
is Aronsson's equation $-\Delta_\infty u =0$ is not surprising.
Infinity harmonic functions (solutions to $-\Delta_\infty u =0$)
appear naturally as limits of $p-$harmonic functions (solutions to
$\Delta_p u = \mbox{div} (|\nabla u|^{p-2} \nabla u) =0$),
\cite{BBM}, and have applications to optimal transport problems,
\cite{EG}, \cite{GAMPR}, image processing, etc, see the survey
\cite{ACJ}.
\end{remark}
In view of Theorem \ref{theorem2 - uniqueness}, it becomes natural
to inquire what happens if condition \eqref{Comp. Condition - thm
2} is violated. In this direction, we show that \eqref{Comp.
Condition - thm 2} is a necessary and sufficient condition for
uniqueness to problem \eqref{limit.problem}. Indeed, if
\eqref{Comp. Condition - thm 2} does not hold, we manage to find
multiple solutions for problem \eqref{limit.problem}.
Nevertheless, we could prove the existence of a minimal one.
\begin{theorem}\label{theorem3 - Multiplicity} Assume that
\eqref{Comp. Condition - thm 2} does not
hold,
then there
exists infinitely many minimizers for the limit problem
\eqref{limit.problem}.
The function
$$
u_\infty (x) =  \max\limits_{y \in \partial
\Omega}
\left (f(y)
        - \Lip (f) |x-y| \right )_{+},
$$
is a minimizer with measure of its positivity set
$$
    \{ u_\infty >0\} = \bigcup\limits_{x\in \partial \Omega}
    B_{\frac{f(x)}{\Lip(f)}}(x) \cap \Omega
$$
strictly less than $\alpha$. Moreover, $u_\infty$ is the minimal
solution, in the sense that any minimizer $v_\infty$ verifies $
v_\infty (x) \geq u_\infty (x)$.
\end{theorem}

\begin{remark} Note the support of the minimal minimizer for problem
\eqref{limit.problem} is given by the set $
 \Omega^\star := \bigcup\limits_{x\in \partial \Omega}
    B_{\frac{f(x)}{\Lip (f) }}(x) \cap \Omega.
$
\end{remark}

Finally, we study geometric properties of the limiting free
boundary, $\partial \{u_\infty > 0 \}$, as well as convergence
issues of the free boundaries $\partial \{u_p > 0 \}$. The next
Theorem we state shows that the the limiting free boundary enjoys
the appropriate geometric features suitable for the study of its
geometric measure properties.

\begin{theorem}\label{Lip and Nondeg of Limit}
    Let $u_p$ be extremals to problem \eqref{problem.p} and assume $u_p \to u_\infty$.
    Then $u_\infty$ is uniformly Lipschitz continuous in $\Omega$, growths linearly away from
    the free boundary and is strongly nondegenerate. That is, for
    a constant $\gamma > 0$,
    $$
        u_\infty(x) \ge \gamma \, \dist (x, \partial \{u_x > 0 \} ), \quad
        \forall x \in \Omega_\infty := \{u_\infty > 0 \}
    $$
    and for any fixed free boundary point $x_0 \in \partial \{ u_\infty > 0 \}$,
    there holds
    $$
        \sup\limits_{B_r(x_0)} u_\infty \ge \gamma r.
    $$
\end{theorem}

The strategy for showing Theorem \ref{Lip and Nondeg of Limit} is
to revisit the $p$-Dirichlet optimization problem
\eqref{problem.p} and verify that these properties hold uniformly
in $p$. As a byproduct of this analysis, we obtain convergence of
the free boundaries $\partial \{u_p > 0\}$ in the Hausdorff
metric.

\begin{theorem}\label{Hausdorff Convergence.intro}
Let $u_p$ be a sequence of minimizers for problem
\eqref{problem.p} and assume $u_p \to u_\infty$, solution to
\eqref{limit.problem}. Then
$$
    \partial \{ u_p > 0 \} \longrightarrow \partial \{ u_\infty >
    0 \}, \qquad  \mbox{ as } p \to \infty,
$$
in the Hausdorff distance.
\end{theorem}

The variational optimization problem \eqref{problem.p} relates, to
some extent, to Bernoulli-type problems governed by the
$p$-Laplacian operator. This is done through a constant free
boundary condition proven to hold for minimizers of problem
\eqref{problem.p}. Indeed, it has been shown (see \cite{FMW06},
\cite{OT}) that $|\nabla u_p| = \lambda_{u_p}$ for a positive
constant $\lambda_{u_p}$ along its free boundary $\partial \{u_p >
0 \}$. This is the so called free boundary condition for the
optimization problem \eqref{problem.p}: a key information when
studying geometric measure as well as smoothness properties of the
free boundary. In this direction we have proven the following
convergence of free boundary conditions.

\begin{theorem}\label{FBCond} Let $u_p$ be a sequence of minimizers for problem
\eqref{problem.p} and  $|\nabla u_p| = \lambda_p$ along $\partial
\{u_p > 0 \}$. Denote $\Omega_\infty := \{ u_\infty > 0 \}$. Then,
up to a subsequence, $(u_p, \lambda_p) \to (u_\infty,
\lambda_\infty)$, with $0 < \lambda_\infty < \infty$ and
$$
    \lim\limits_{\begin{array}{l} x\to \partial \Omega_\infty \\ x \in \Omega_\infty
    \end{array}}
    \dfrac{u_\infty(x)}{\dist (x,\partial \Omega_\infty)} =
    \lambda_\infty.
$$
\end{theorem}

When $\Omega$ is convex and $f$ is constant, Theorem \ref{FBCond}
can be seen in connection to the results of Manfredi, Petrosyan
and Shahgholian, \cite{MPS}, who study convergence issues, as $p
\to \infty$, for Bernoulli-type problems.

\medskip

The rest of the paper is organized as follows: in the next section
we prove Theorem \ref{theorem 1 - convergence}; in Section~3 we
study the limit problem under condition \eqref{Comp. Condition -
thm 2} and in Section~4 we deal with the complementary case;
finally in Section~5 we include some uniform bounds for the
sequence $u_p$ (showing uniform non-degeneracy of the free
boundary) and we study the convergence of the free boundaries.

\section{Proof of Theorem \ref{theorem 1 - convergence}}

In this section we prove Theorem \ref{theorem 1 - convergence}.
The main issue of the proof is to find bounds for the energy
$\left(\int_\Omega |\nabla u_p |^p \right)^{1/p}$ of a minimizer
that are independent of $p$.

\medskip

\noindent \textbf{Proof of Theorem \ref{theorem 1 - convergence}}
Let us fix hereafter a Lipschitz extension of $f$, which we shall
denoted by $v$, among functions within
\begin{equation}\label{K.infty}
 K_\infty = \left\{ \varphi \in W^{1,\infty} (\Omega) \suchthat \ \varphi =f, \ \mbox{on} \
 \partial \Omega, \ |\{\varphi > 0 \}| = \alpha \right\}.
\end{equation}
Clearly, since $\Omega$ is bounded,  $v$ competes in the
minimization problem \eqref{problem.p}. Thus using $v$ as a test
function in problem \eqref{problem.p} we obtain
$$
    \left(\int_\Omega |\nabla u_p |^p \right)^{1/p}\leq \left(
    \int_\Omega |\nabla v |^p \right)^{1/p} \leq \Lip(v)
    |\Omega|^{1/p} \leq C,
$$
where $C$ is a constant independent of $p$. With exponent $q<
\infty$ fixed, we now argue as follows
$$
    \left(\int_\Omega |\nabla u_p |^q \right)^{1/q}\leq
    \left(\int_\Omega |\nabla u_p |^p \right)^{1/p}
    |\Omega|^{p/(q(p-q))} \leq \Lip(v) |\Omega|^{1/p + p/(q(p-q))} \leq
    C.
$$
Therefore, the sequence $u_p$ is uniformly bounded in $W^{1,q}
(\Omega)$, and its weak limit as $p\to \infty$,  $u_\infty$
verifies
$$
\left(\int_\Omega |\nabla u_\infty |^q \right)^{1/q}\leq \Lip(v)
|\Omega|^{1/q} \leq C.
$$
Taking $q\to \infty$ and performing a diagonal argument,  we
obtain a subsequence $u_p$ that converges weakly in every $W^{1,q}
(\Omega)$, $1<q<\infty$ to a limit $u_\infty \in W^{1,\infty}
(\Omega)$ such that
$$
    \| \nabla u_\infty \|_{L^\infty (\Omega)} \leq \Lip(v).
$$
Let us now turn our attention towards estimating the Lebesgue
measure of $\{ u_\infty > 0 \}$. Fixed an $\epsilon > 0$, thanks
to the uniform convergence, for $p$ large enough, there holds
$$
    \left \{ u_\infty > \epsilon  \right \} \subset \{ u_p > 0 \}.
$$
Hence we conclude that
$$
\Leb \left ( \{ u_\infty >0 \} \right ) = \lim_{\epsilon \to 0}
\Leb \left ( \{ u_\infty > \epsilon \} \right ) \leq \alpha.
$$
Therefore, we have proved that $u_\infty$ is an extremal for the
limit problem \eqref{limit.problem}.

It remains to prove that $u_\infty $ is indeed $\infty$-harmonic
in its set of positivity. Following \cite{CIL} let us recall  the
definition of viscosity solution.

\begin{defi} \label{def.sol.viscosa}
Consider the boundary value problem
\begin{equation}\label{ec.viscosa.con.borde}
\begin{array}{ll}
F (x,D u, D^2 u )   = 0 \qquad & \mbox{in } \Omega.
\end{array}
\end{equation}
\begin{enumerate}
\item A lower semi-continuous function $ u $ is a viscosity supersolution if for every
$ \phi \in C^2(\overline{\Omega})$ such that $ u-\phi $ has a
strict minimum at the point $ x_0 \in \Omega$ with $u(x_0)=
\phi(x_0)$ we have:
$$
F(x_0, D \phi (x_0), D^2\phi (x_0)) \ge 0.
$$
\item An upper semi-continuous function $u$ is a subsolution if for every $ \phi  \in
C^2(\overline{\Omega})$ such that $ u-\phi $ has a strict maximum
at the point $ x_0 \in \Omega$ with $u(x_0)= \phi(x_0)$ we have:
$$
F(x_0, D \phi (x_0), D^2\phi (x_0)) \le 0 .
$$
\item Finally, $u$ is a viscosity solution if it is a super and a
subsolution.
\end{enumerate}
\end{defi}

If we have a weak $p-$harmonic function (in the sense of
distribution) that is continuous then it is a viscosity solution.
This is the content of our next result.

\begin{lema} \label{lema.sol.debil.es.sol.viscosa}
Let $u$ be a continuous weak solution of $\Delta_p u =0$ in some
domain $\Omega$ for $p>2$. Then $u$ is a viscosity solution of
\begin{equation}\label{ec.viscosa.con.borde.p}
\begin{array}{ll}
-(p-2) |D u|^{p-4}\Delta_\infty u - |D u|^{p-2}
\Delta u = 0
 \qquad & \mbox{in } \Omega.\\
\end{array}
\end{equation}
\end{lema}

\begin{proof}  Let $x_0 \in \Omega$ and a test function $\phi$ such
that $u(x_0)=\phi (x_0)$ and $u-\phi$ has a strict minimum at
$x_0$. We want to show that
$$
-(p-2) |D \phi|^{p-4}\Delta_\infty \phi (x_0) - |D \phi|^{p-2}
\Delta \phi (x_0) \ge 0.
$$
Assume that this is not the case, then there exists a radius $r>0$
such that
$$
-(p-2) |D \phi|^{p-4}\Delta_\infty \phi (x) - |D \phi|^{p-2}
\Delta \phi (x) < 0,
$$
for every $x\in B(x_0,r)$. Set $m = \inf_{|x-x_0|=r} (u-\phi)(x)$
and let $\psi (x) = \phi(x) + m/2$. This function $\psi$ verifies
$\psi (x_0) > u(x_0)$ and
$$
-\mbox{div} ( |D \psi|^{p-2} D \psi) <0.
$$
Multiplying by $(\psi - u)^+$ extended by zero outside $B(x_0,r)$
we get
$$
\int_{\{ \psi > u \}} |D \psi|^{p-2} D \psi D (\psi - u) < 0.
$$
Taking $(\psi - u)^+$ as test function in the weak form of the
equation we get
$$
\int_{\{ \psi > u \}} |D u|^{p-2} D u D (\psi - u) =0.
$$
Hence,
$$
\displaystyle C(N,p)\int_{\{ \psi > u \}} |D \psi - D u|^{p}  \le \int_{\{ \psi
> u \}} \langle |D \psi|^{p-2} D \psi - |D u|^{p-2}
D u, D (\psi - u) \rangle <0,
$$
a contradiction. This proves that $u$ is a viscosity
supersolution. The proof of the fact that $u$ is a viscosity
subsolution runs as above, we omit the details.
\end{proof}

We are now  ready to prove that the limit $\lim_{p_i\to \infty}
u_{p_i} = u_{\infty}$ satisfies the desired PDE in its set of
positivity. In fact, let us check that $-\Delta_\infty
u_{\infty}=0$ in the viscosity sense in the set $\{u_\infty >0
\}$.  Let us recall the standard proof. Let $\phi$ be a smooth
test function such that $u_{\infty}-\phi$ has a strict maximum at
$x_0 \in \{u_\infty >0 \}$. Since $u_{p_i}$ converges uniformly to
$u_{\infty}$ we get that $u_{p_i} - \phi$ has a maximum at some
point $x_i\in \Omega$ with $x_i \to x_0$ and moreover we have that
$u_{p_i} >0$ in a whole fixed neighborhood of $x_0$ (and therefore
$u_{p_i} (x_i )
>0$ and every $u_{p_i}$ is $p-$harmonic there). Next, we use the
fact that $u_{p_i}$ is a viscosity solution of $ -\Delta_p u_{p}
=0 $ in the set $\{u_{p_i} >0 \}$ and we obtain
\begin{equation}\label{sol.viscosa.p.2}
-(p_i-2) |D \phi|^{p_i-4} \Delta_\infty \phi (x_i) - |D
\phi|^{p_i-2} \Delta \phi (x_i) \le 0.
\end{equation}
If $D \phi (x_0)=0$ we get $-\Delta_\infty \phi (x_0) \le 0$. If
this is not the case, we have that $D \phi (x_i) \neq 0$ for large
$i$ and then
$$
- \Delta_\infty \phi (x_i) \le  \frac{1}{p_i-2} |D \phi|^{2}
\Delta \phi (x_i) \to 0, \quad \mbox{as } i \to \infty .
$$
We conclude that
$$
-\Delta_{\infty} \phi (x_0) \le 0.
$$
That is $u_{\infty}$ is a viscosity subsolution of
$-\Delta_{\infty} u_{\infty} =0$.

A similar argument shows that $u_{\infty}$ is also a supersolution
and therefore a solution of $-\Delta_\infty u_{\infty} =0$ in
$\Omega$. The proof of Theorem \ref{theorem 1 - convergence} is
completed. \hfill $\square$

\section{Proof of Theorem \ref{theorem2 - uniqueness}}

In this section we deal with the situation in which we have
uniqueness for the limit problem. We will assume that condition
\eqref{Comp. Condition - thm 2} holds, that is,
$$
    \Leb \left( \bigcup_{y \in\partial \Omega} B_{\frac{f(y)}{\Lip (f)}}
    (y) \cap \Omega \right) \geq \alpha.
$$

Note that with the notations of the statement of Theorem
\ref{theorem2 - uniqueness} this implies that
$$
\lambda^{\star} \geq \Lip (f).
$$
This fact is crucial in the course of next proof.

\medskip

\noindent \textbf{Proof of Theorem \ref{theorem2 - uniqueness}}
     Let $v_\infty$ be a minimizer for problem
     \eqref{limit.problem}. Existence of such a minimizer is
     assured by Theorem \ref{theorem 1 - convergence}. Let us denote
     $$
        \Omega_\infty := \{ v_\infty > 0 \} \subset \Omega.
     $$
     For each free boundary point $y \in \partial \Omega_\infty$,
     let $x\in \partial \Omega$ be a point satisfying
     $$
        |x - y| = \dist (y, \partial \Omega).
     $$
     Using the Lipschitz continuity of $v_\infty$, we obtain the
     following estimate
     \begin{equation}\label{eq01 - lemma2}
        f(x) \le \Lip(v_\infty) |x-y|.
     \end{equation}
    From \eqref{eq01 - lemma2}, we immediately conclude that
    \begin{equation}\label{eq02 - lemma2}
        \bigcup\limits_{x\in \partial \Omega}
        B_{\frac{f(x)}{\Lip(v_\infty)}}(x) \cap \Omega \subset \Omega_\infty,
    \end{equation}
    which, in particular, implies
    \begin{equation}\label{eq03 - lemma2}
        \Leb \left ( \bigcup\limits_{x\in \partial \Omega}
        B_{\frac{f(x)}{\Lip(v_\infty)}}(x) \cap \Omega \right ) \le \alpha =
        \Leb \left ( \bigcup\limits_{x\in \partial \Omega}
        B_{\frac{f(x)}{\lambda^{\star}}}(x) \cap \Omega \right ).
    \end{equation}
    From above we obtain
    \begin{equation}\label{eq04 - lemma2}
        \lambda^{\star} \le \Lip(v_\infty).
    \end{equation}

On the other hand, let
$$
    \Omega^\star := \bigcup\limits_{x\in \partial \Omega}
    B_{\frac{f(x)}{\lambda^{\star}}}(x) \cap \Omega.
$$
Then, $u_\infty$, defined as the solution to
$$
    \left \{
        \begin{array}{rll}
            \Delta_{\infty} u_\infty &=& 0 \qquad \text{ in } \Omega^\star, \\
            u_\infty & = & f \qquad \text{ on } \partial \Omega,  \\
            u_\infty & = & 0 \qquad \text{ on } \partial \Omega^\star
            \cap \Omega
        \end{array}
    \right.
$$
competes in the minimization problem \eqref{limit.problem}, thus
    \begin{equation}\label{eq05 - lemma2}
        \Lip(u_\infty) \ge \Lip(v_\infty)
    \end{equation}
    In the sequel, we will use the fact that $u_\infty$ is the best Lipschitz
    extension of the boundary data $f$ on $\partial \Omega$ and
    $0$ on $\partial \Omega^\star \cap \Omega$ together with the geometric
    compatibility condition \eqref{Comp. Condition - thm 2} to
    bridge these inequalities. For that we consider the auxiliary barrier function
    $$
        \psi(x) := \max\limits_{y \in \partial \Omega} \left (f(y)
        - \lambda^{\star} |x-y| \right )_{+}.
    $$

    We initially verify that
    \begin{equation}\label{eq06 - lemma2}
        \Lip(\psi) = \lambda^{\star}.
    \end{equation}
    Indeed, let $x_1$ and $x_2$ be two points in
    $\Omega$. We assume $0 < \psi(x_1) < \psi(x_2)$. Let $y_1$ and
    $y_2$ be such that
    $$
        \psi(x_i) = f(y_i) - \lambda^{\star} |x_i-y_i|, \quad \text{i=1,2}.
    $$
    From the definition of $\psi$, we know
    $$
        \psi(x_1) \ge f(y_2) - \lambda^{\star} |x_1-y_2|.
    $$
    We now estimate
    $$
        \begin{array}{lll}
            0 < \psi(x_2) - \psi(x_1) &\le& f(y_2) - \lambda^{\star} |x_2-y_2|
            - \left (f(y_2) - \lambda^{\star} |x_1-y_2| \right ) \\
            &\le& \lambda^{\star} \left ( |x_1-y_2| - |x_2-y_2| \right ) \\
            & \le & \lambda^{\star} |x_1 - x_2|.
        \end{array}
    $$

    Our next step is to check that $\psi$ matches the desired
    boundary conditions. Well, it is clear from its definition that
    $$
        \psi \big |_{\partial \Omega^\star} = 0.
    $$
    Proving $\psi$ agrees with $f$ on $\partial \Omega$ is
    equivalent to showing that
    \begin{equation}\label{eq07 - lemma2}
        f(x) = \max_{y \in \partial \Omega} \left\{ (f(y) - \lambda^{\star} |x-y|)_+
        \right\}, \quad \forall x \in \partial \Omega.
    \end{equation}
    Let us assume, for sake of contradiction, that \eqref{eq07 - lemma2} does not hold.
    This would readily imply that there exist two points $x$, $y$ on $\partial \Omega$ with
    $$
        \lambda^{\star} |x-y| < f(y) - f(x).
    $$
    That is,
    $$
        \lambda^{\star} < \Lip(f) = \sup_{x,y \in \partial \Omega} \left\{ \frac{
        |f(x) - f(y)|}{|x-y|} \right\}.
    $$
    which contradicts (\ref{Comp. Condition - thm 2}).

As a remark, note that when we take two pints $x,y \in \partial
\Omega$ we get
$$
|\psi (x) - \psi (y)| = |f(x) - f(y)| \leq \Lip (f) |x-y|
$$
Thus $\Lip(\psi) = \max(\lambda^\star, \Lip(f)) =
\lambda^{\star}$.

    Once verified that $\psi$ has the same boundary condition as $u_\infty$,
    from the fact that $u_\infty$ is the best Lipschitz extension its boundary data, we
    know
    \begin{equation}\label{eq08 - lemma2}
        \Lip (u_\infty) \le \Lip(\psi) = \lambda^{\star}.
    \end{equation}

Now let us show that $u_\infty$ coincides with the barrier
$$
\psi(x) = \max\limits_{y \in \partial \Omega} \left (f(y)
        - \lambda^{\star} |x-y| \right )_{+}.
$$
We have that $ u_\infty$ is a minimizer for the limit problem,
hence we must have
$$
u_\infty (x) \geq  \max\limits_{y \in \partial \Omega} \left (f(y)
        - \lambda^{\star} |x-y| \right )_{+}.
$$
In fact, assume that this is not the case, then there exists $x_0$
such that $u_\infty (x_0) < \psi (x_0)$. Now, considering
quotients that involve $x_0$ and points on $\partial
\Omega$ we can easily conclude that $\Lip (u_\infty)
>
\lambda^{\star} = \Lip (\psi)$, a contradiction since $\psi$ is a
competitor in the limit problem.

Therefore, we obtain that both functions have the same positivity
set (both sets have the same measure and one is included in the
other).

Now, arguing as before, assume that there exists $x_0$ such that $
u_\infty (x_0) > \psi (x_0) $. In this case, comparing quotients
defining the Lipschitz constant with $x_0$ and points on the
boundary of the positivity set, we get $ \Lip (u_\infty ) >
\lambda^{\star} = \Lip (\psi) $. This contradicts again the fact
that $u_\infty$ is optimal for the limit problem.

Combining (\ref{eq02 - lemma2}), (\ref{eq03 - lemma2}) (\ref{eq04
- lemma2}) and (\ref{eq08 - lemma2}), together with the fact that
$u_\infty$ and $\psi$ are $\infty$-harmonic in $\Omega^*$ with the
same value on the boundary of this set, we finish up the proof of
Theorem
\ref{theorem2 - uniqueness}. \hfill $\square$

\section{Proof of Theorem \ref{theorem3 - Multiplicity}}

Now let us show that when (\ref{Comp. Condition - thm 2}) does not
hold there is no uniqueness for minimizers of the limit problem.

\medskip

\noindent \textbf{Proof of Theorem \ref{theorem3 - Multiplicity}}
As before, let $\lambda^{\star}$ be such that
$$
    \Omega^\star := \bigcup\limits_{x\in \partial \Omega}
    B_{\frac{f(x)}{\lambda^{\star}}}(x) \cap \Omega
$$
has Lebesgue measure precisely $\alpha$ and assume that
(\ref{Comp. Condition - thm 2}) does not hold, that is,
$$
\Lip (f) > \lambda^{\star}.
$$
Let
$$
D := \bigcup\limits_{x\in \partial \Omega}
    B_{\frac{f(x)}{\Lip (f) }}(x) \cap \Omega.
$$
We have
$$
\Leb (D) < \alpha.
$$
By our previous result we have that
$$
        \psi(x) := \max\limits_{y \in \partial \Omega} \left (f(y)
        - \Lip (f) |x-y| \right )_{+}
$$
is an extremal for the limit problem with measure $\Leb (D)$.

Now, let $v_\infty$ be an extremal for the limit problem with
measure $\alpha$. Then, as $v_\infty = f$ on $\partial \Omega$ we
have
$$
    \Lip (v_\infty ) \geq \Lip (f) = \Lip (\psi).
$$

On the other hand $\psi$ is a competitor in the limit problem with
measure $\alpha$ and hence
$$
    \Lip (\psi ) \geq \Lip (v_\infty ).
$$
We conclude that
$$
    \Lip (\psi) = \Lip (v_\infty ) = \Lip (f)
$$
and then $\psi$ is also a maximizer for the limit problem.

Moreover, we have that
$$
    \psi (x) \leq v_\infty (x), \qquad x \in D,
$$
if not the Lipschitz constant of $v_\infty$ is greater than
$\Lip(\psi)$. Indeed, let us assume that there exists $x_0 \in D$
such that
$$
\psi (x_0) > v_\infty (x_0).
$$
That is,
$$
\max\limits_{y \in \partial \Omega} \left (f(y)
        - \Lip (f) |x_0-y| \right )_{+} > v_\infty (x_0)
$$
From where we get that there exists $y \in \partial \Omega$ such
that
$$
    f(y) - \Lip (f) |x_0-y| > v_\infty (x_0)
$$
that is to say that (using that $v_\infty = f$ on $\partial
\Omega$),
$$
v_\infty (y) - v_\infty (x_0) > \Lip (f) |x_0-y|
$$
that clearly implies that
$$
    \Lip (f) < \Lip (v_\infty).
$$

Therefore, we have that $\psi$ is the {\it minimal} extremal for
the limit problem and hence we obtain the following estimate for
the support of any extremal $v_\infty$,
$$
D := \bigcup\limits_{x\in \partial \Omega}
    B_{\frac{f(x)}{\Lip (f) }}(x) \cap \Omega
\subset \{ v_\infty >0\}.
$$

Now, let
$$
D_\delta := \bigcup\limits_{x\in \partial \Omega}
    B_{\frac{f(x)}{\Lip (f) }}(x) \cap \Omega
+ B(0,\delta),
$$
with $\delta$ small such that
$$
\Leb (D_\delta) < \alpha.
$$
In this set $D_\delta$, let us consider $v_\infty$ the solution to
$$
    \left \{
        \begin{array}{rll}
            \Delta_{\infty} v_\infty &=& 0 \qquad \text{ in } D_\delta, \\
            v_\infty & = & f \qquad \text{ on } \partial \Omega,  \\
            v_\infty & = & 0 \qquad \text{ on } \partial D_\delta.
        \end{array}
    \right.
$$
Since $D\subset D_\delta$, we have
$$
    \Lip (v_\infty) = \Lip (f).
$$
To prove this fact, let us consider in the set $D_\delta$ the
boundary value
$$
F(x) = \left\{ \begin{array}{ll} f(x) \qquad & x \in \partial
\Omega, \\
0 \qquad & x \in \partial D_\delta \cap \Omega.
\end{array}\right.
$$
This boundary datum $F$ is a Lipschitz function with Lipschitz
constant given by
$$
    \Lip (F) = \sup_{x,y \in \partial D_\delta} \frac{|F(x) -
F(y)|}{|x-y|}.
$$
Let us estimate this Lipschitz constant $\Lip (F)$. If $x,y \in
\partial D_\delta \cap \Omega$ then
$$
\frac{|F(x) - F(y)|}{|x-y|} =0 < \Lip(f).
$$
When $x,y \in \partial \Omega$, clearly
$$
\frac{|F(x) - F(y)|}{|x-y|} \leq \Lip(f).
$$
And finally when $x\in \partial \Omega$ and $y \in \partial
D_\delta \cap \Omega$ we have
$$
\frac{|F(x) - F(y)|}{|x-y|}  = \frac{|f(x)}{|x-y|}  < \Lip(f).
$$
We are using the fact that $D \subset \Omega_\delta$ and hence the
distance $|x-y|$ is bigger than $f(x)/\Lip (f)$, to see this fact,
just take $y\in \partial D$ then for any $x \in \partial \Omega$
we have
$$
f(x) - \Lip (f)  |x-y| \leq 0
$$
that is to say
$$
|x-y| \leq \frac{f(x)}{ \Lip (f)}.
$$
Therefore we conclude that
$$
    \Lip (F) = \Lip (f),
$$
and since $v_\infty$ has the same Lipschitz constant as $F$ (it is
its best possible Lipschitz extension) we conclude that
$$
    \Lip (v_\infty) = \Lip (f).
$$

Hence $v_\infty$ is also an extremal for the limit problem that is
positive on $\partial D \subset (D_\delta)^o$ (the strong maximum
principle holds for $\infty$-harmonic functions) and hence we
conclude that $v_\infty \neq \psi$.

With these estimates we can conclude that there is no strict
monotonicity with respect to the measure in the limit problem.
\hfill $\square$

\medskip

Now, we can state further consequences of our previous results.

\begin{theorem}\label{theoremPepe-section4}
Assume that
$$
\beta := \Leb \left( \bigcup\limits_{x\in \partial \Omega}
    B_{\frac{f(x)}{\Lip (f) }}(x) \cap \Omega \right) < \alpha.
$$
Then we have
$$
\lim_{p\to \infty} \mathfrak{P}_p (\alpha) = \lim_{p\to \infty}
\mathfrak{P}_p (\beta)
$$
in the sense that if $u_p$ is an extremal for $\mathfrak{P}_p
(\alpha)$ and $v_p$ is an extremal for $\mathfrak{P}_p (\beta)$
then
$$
\lim_{p \to \infty} \left( \int_{\Omega} |\nabla u_p (X)|^p dX
\right)^{1/p} = \lim_{p \to \infty} \left( \int_{\Omega} |\nabla
v_p (X)|^p dX \right)^{1/p}
$$
and moreover,
$$
v_p \to \psi \qquad \mbox{ and } \qquad u_p
\to  u_\infty
$$
uniformly in $\overline{\Omega}$ with
$$
    \Lip (u_\infty) = \Lip (\psi) = \Lip (f)\qquad \mbox{ and } \qquad
\psi (x) \leq u_\infty (x).
$$
\end{theorem}

One possible conclusion of this fact is that the boundary datum
$f$ is so that the limit problem has many solutions and hence we
are ``wasting measure" when considering the problem with $\alpha$
instead of $\beta$. In fact, the value of the minimum for
$\mathfrak{P}_p (\alpha)$ and for $\mathfrak{P}_p (\beta)$ are
almost the same for $p$ large and the minimal solution of the
limit problem is $\psi$ (which is the unique minimizer for
$\mathfrak{P}_\infty (\beta)$).

\section{Uniform estimates and free boundary convergence issues}

This section is devoted to establish Theorems \ref{Lip and Nondeg
of Limit}, \ref{Hausdorff Convergence.intro} and \ref{FBCond}. For
that we shall revisit the study of the $p$-Dirichlet energy
minimization problem with volume constraint, \eqref{problem.p}
carried out in \cite{OT} and in \cite{FMW06}. Our strategy is to
seize uniform-in-$p$ properties and afterwards explore their
impact on the limiting problem \eqref{limit.problem}.
\par
It is well established in the literature that ordinary techniques
from the Calculus of Variations are not suitable to approach
directly optimal design problems with volume constraints. Indeed,
to establish existence of a minimizer for Problem
\eqref{problem.p} requires a careful analysis, involving penalty
method and geometric measure perturbation techniques.
\par
Penalization version of problem \eqref{problem.p} can be easily
set-up. Indeed, for each $L > 0$, let
\begin{equation}\label{penalization term}
    \varrho_L(t) := L (t - \alpha)^{+}.
\end{equation}
We then define the $L$-penalized problem for the $p$-Dirichlet
integral, as
\begin{equation}\label{penalized problem}\tag{$\mathfrak{P}_p^L$}
    \min \left \{\int_{\Omega} |\nabla u (X)|^p dX + \varrho_L\left (\{u>0\} \right )
    \suchthat u \in W^{1,p} (\Omega), \ u = f \mbox{ on } \partial \Omega
\right\}.
\end{equation}
Notice that problem \eqref{penalized problem} does not involve
volume constraint anymore, thus the proof of existence of a
minimizer, $u_p^L$, for problem \eqref{penalized problem} follows
a standard scheme from the Calculus of Variations. It is also
simply to check that $u_p^L \ge 0$ and $\Delta_p u_p^L$ is a
non-negative Radon measure supported on $\partial \{ u_p^L > 0
\}$. In particular, $u_p^L$ is $p$-harmonic in its set of
positivity, that is, $u_p^L$ satisfies the following PDE
$$
    \Delta_p u_p^L = 0, \qquad  \text{ in } \{u_p^L > 0\}.
$$
Although locally $C^{1,\alpha}$ within $\{u_p^L > 0 \}$, notice
that Lipschitz is the optimal regularity for $u_p^L$ in $\Omega$.
This is because $\nabla u_p^L$ jumps from positive slope to zero
along the free boundary $\partial \{u_p^L > 0\}$. Indeed it has
been proven in \cite{OT, FMW06} that for each $L$ fixed $u_p^L$ is
locally Lipschitz continuous in $\Omega$. Our next lemma gives the
precise dependence of the Lipschitz norm of $u_p^L$ with respect
to $p$ and the penalty charge $L$. This lemma is essentially taken
from \cite{Teix-Preprint}. We present a proof here as a courtesy
to the readers.
\begin{lemma}\label{Lip contiunity} Let $u_p^L$ be a minimizer for
\eqref{penalized problem}. Then,
$$
    \| \nabla u_p^L \|_{L^\infty(\Omega)} \le C L^{1/p},
$$
where $C$ is a constant that depends only on dimension, $f$ and
$\alpha$.
\end{lemma}

\begin{proof}

Since we are interested in the limiting problem, we will only deal
with the case $p \gg 1$. We will follows the approach suggested in
\cite{AC}, keeping track of the precise constants that appear on
the estimates. From the minimality of $u_p^L$, we deduce, for any
ball $B = B_d(x_0) \subset \Omega$, centered at a free boundary
point, i.e., $x_0 \in
\partial \{u_p^L
> 0 \}$, there holds
\begin{equation} \label{Lip Eq01}
     L \cdot \Leb \left ( \{ x \in B_d(x_0) \suchthat u_p^L(x) = 0 \} \right ) \ge
     c_0  \left ( \int_{\Omega} \left |
    \nabla \left (u_p^L - \mathfrak{h}_p \right )(x) \right|^p  dx \right ),
\end{equation}
where $\mathfrak{h}_p$ is the $p$-harmonic function in $B_d(x_0)$
that agrees with $u_p^L$ on $\partial B_d(x_0)$ and $c_0$ is a
constant that depends only upon dimension. For any direction
$\nu$, we define
$$
    r_\nu := \min \left \{ r \suchthat \frac{1}{4} \le r \le 1
    \textrm{ and } u_p^L(x_0 + d r \nu) = 0 \right \}
$$
if such a set is nonempty; otherwise, we put $r_\nu = 1$. Taking
into account that $$u_p^L(x_0 + d r_\nu \nu) = 0$$ whenever $r_\nu
< 1$, we can compute,
\begin{equation}\label{Lip Eq02}
    \begin{array}{lll}
        \mathfrak{h}_{p}(x_0 + d r_\nu \nu) & = & \displaystyle \int_{r_\nu}^1 \dfrac{d}{dr}
        (u_p^L - \mathfrak{h}_{p})(x_0 + d r \nu) dr \\
        & \le & d \cdot (1-r_\nu)^{1/p'} \times \left [ \displaystyle \int_{r_\nu}^1
        |\nabla ( \mathfrak{h}_{p} - u_p^L)(x_0 +  r\nu)|^p dr \right
        ]^{1/p}.
    \end{array}
\end{equation}
Here $\frac{1}{p} + \frac{1}{p'} = 1$. Now, by the Harnack
inequality, we know
\begin{equation}\label{Lip Cont Eq02}
    \inf\limits_{B_{\frac{2}{3}d}(x_0)} \mathfrak{h}_{p} \ge c_1
    \mathfrak{h}_{p}(x_0),
\end{equation}
for a constant $c_1>0$ that depends only on dimension (see, for
instance, \cite{KMV}). Let us consider the following barrier
function, $\mathrm{b}$, given by
\begin{equation}\label{Lip Cont Eq03}
    \left \{
        \begin{array}{rll}
            \Delta_p \mathrm{b} &=& 0\qquad  \textrm{ in } B_1(0) \setminus
            B_{\frac{2}{3}}(0), \\
            \mathrm{b} &=& 0 \qquad \textrm{ on } \partial B_1(0), \\
            \mathrm{b} &=& c_1 \qquad \textrm{ in }
            \overline{B_{\frac{2}{3}}(0)},
        \end{array}
    \right.
\end{equation}
where $c_1$ is the universal constant in (\ref{Lip Cont Eq02}). By
the Hopf's maximum principle, there exists a universal constant
$c_2 >0$, depending only on dimension, such that
\begin{equation}\label{Lip Cont Eq04}
    b(x) \ge c_2 \left (1 - |x| \right ).
\end{equation}
By the maximum principle and (\ref{Lip Cont Eq04}) we can write
\begin{equation}\label{Lip Cont Eq05}
    \mathfrak{h}_{p}(x_0 + dx) \ge \mathfrak{h}_{p}(x_0)\cdot
    b(x) \ge c_2 \mathfrak{h}_{p}(x_0)\cdot  (1 - |x|).
\end{equation}
Combining (\ref{Lip Eq01}) and (\ref{Lip Cont Eq05}) we end up
with
\begin{equation}\label{Lip Cont Eq06}
    d^p \cdot
    \left [ \displaystyle \int_{r_\nu}^1|\nabla ( \mathfrak{h}_{p} - u_p^L )
    (x_0 +  r\nu)|^p dr \right ] \ge c_3 \mathfrak{h}_{p}^p(x_0) \cdot  (1 -
    r_\nu).
\end{equation}
Integrating (\ref{Lip Cont Eq06}) with respect to $\nu$ over
$\mathbb{S}^{n-1}$, taking into account the definition of $r_\nu$,
we find
\begin{equation}\label{Lip Cont Eq07}
    \left ( \dfrac{\mathfrak{h}_{p}(x)}{d} \right )^p
    \cdot \int_{B_d(x) \setminus B_{d/4}(x)} \chi_{\{ u_p^L =
    0\}} dx \le C_4 \int_{B_d(x)} \left |\nabla \left (\mathfrak{h}_{p} -  u_p^L
    \right )(x) \right |^p dx.
\end{equation}
If we replace, in all of our arguments so far, $B_{d/4}(x)$ by
$B_{d/4}(\overline{x})$, for any $\overline{x} \in \partial
B_{d/2}(x)$, we obtain
\begin{equation}\label{Lip Cont Eq08}
    \left ( \dfrac{\mathfrak{h}_{p}(x)}{d} \right )^p
    \cdot \int_{B_d(x) \setminus B_{d/4}(\overline{x})} \chi_{\{ u_p^L =
    0\}} dx \le \tilde{C}_4 \int_{B_d(x)} \left |\nabla \left (\mathfrak{h}_{p} -  u_p^L
    \right )(x) \right |^p dx,
\end{equation}
for every $\overline{x} \in \partial B_{d/2}(x)$.

Integrating (\ref{Lip Cont Eq08}) with respect to $\overline{x}$,
yields:
\begin{equation}\label{Lip Cont Eq09}
    \left ( \dfrac{\mathfrak{h}_{p}(x)}{d} \right )^p
    \cdot \left | \left \{ x \in B_d(x) \suchthat u_p^L (x) = 0  \right \} \right |
    \le C_5 \int_{B_d(x)} \left |\nabla \left (\mathfrak{h}_{p} -  u_p^L
    \right )(x) \right |^p dx.
\end{equation}
Now we argue as follows: let $\rho := \dist(x, \partial \{u_p^L
> 0 \})$ and for each $0 < \delta <\!\!< 1$, denote
$\mathfrak{h}_p^\delta$ the $p$-harmonic function in $B_{\rho +
\delta}(x)$ that agrees with $u_p^L$ on $\partial B_{\rho +
\delta}(x)$. Combining (\ref{Lip Eq01}) and (\ref{Lip Cont Eq09})
together with standard elliptic estimate, we deduce
\begin{equation}\label{Lip Cont Eq10}
        u_p^L(x)  =  \mathfrak{h}_\delta(x) + \text{o}(1)
        \le  C_6L^{1/p} (\rho + \delta) + \text{o}(1), \quad
        \text{as } \delta \searrow 0,
\end{equation}
for a constant $C_6$ that depends on dimension, $f$ and $\alpha$.
Letting $\delta \searrow 0$ in (\ref{Lip Cont Eq10}) we finally
conclude
$$
    u_p^L(x) \le C_6 L^{1/p} \dist \left (x, \partial
    \Omega_\lambda^\star \right ),
$$
which clearly implies that $u_p^L$ is Lipschitz continuous up to
the free boundary $\partial \{ u_p^L > 0 \}$ and $\|\nabla
u_p^L\|_\infty \lesssim L^{1/p}$. Lemma \ref{Lip contiunity} is
proven.
\end{proof}
Another important piece of information concerns uniform
non-degeneracy.
\begin{lemma}\label{nondegeneracy} Let $x \in \{ u_p^L > 0 \}$ be
a free boundary point. Then
    \begin{equation}\label{Nondegeneracy Penalized Problem}
        L^{-1/p} \underline{c} \cdot \dist \left (x, \partial \{u_p^L > 0 \} \right ) \le u_p^L(x)
    \end{equation}
for a constant $\underline{c}$ that depends only upon dimension,
$f$ and $\alpha$. Moreover the following strong non-degeneracy
holds
    \begin{equation}\label{Strong Nondegeneracy Penalized Problem}
        \sup\limits_{B_r(x_0)} u_p^L \ge L^{-1/p} \underline{c}_1 r,
    \end{equation}
for any free boundary point $x_0 \in \partial \{u_p^L > 0 \}$. The
constant $\underline{c}_1$ depends only on dimension, $f$ and
$\alpha$ and is independent of $p$.
\end{lemma}
The proof of Lemma \ref{Nondegeneracy Penalized Problem} is, by
now, classical in variational free boundary theory. It relies on
``cutting" a small hole around the free boundary point and
comparing the result with the original optimal design. For further
details we refer the readers to \cite{Teix-Preprint}, Theorem 6.2.
As observed in the proof of Lemma \ref{Lip contiunity}, the fact
that $\underline{c}$ and $\underline{c}_1$ are universal is a
consequence of uniform-in-$p$ Harnack inequality and
uniform-in-$p$ Hopf boundary maximum principle. We skip the
details here.
\par
The penalty method strategy is based on the idea that if $L$ is
large enough (but still finite), one expects that minimizers for
\eqref{penalized problem} would rather prefer to obey the volume
constraint, $\Leb \left (\{ u_p^L > 0 \} \right ) \le \alpha$.
Therefore it would be a solution for the original problem,
\eqref{problem.p}. Such a strategy does work, \cite{OT, FMW06} and
\cite{Teix-Preprint}, however it relies on a fine geometric
measure perturbation approach. The following theorem is a
consequence of the analysis carried out in \cite{Teix-Preprint},
section 7:
\begin{lemma}\label{est penalty term} There exists a universal constant $C$,
depending only on dimension, $f$ and $\alpha$, but independent of
$p$, such that if
$$
    L \ge Cp,
$$
then
$$
\Leb \left ( \{u_p^{L} > 0\} \right ) \le \alpha.
$$
Therefore, $u_p^{Cp}$ is a solution to problem \eqref{problem.p}.
\end{lemma}

It is important to notice that any minimizer, $u_p$, of problem
\eqref{problem.p} is also a minimizing function for problem
($\mathfrak{P}_p^{Cp}$). As a consequence, combining Lemma
\ref{Lip contiunity}, \ref{nondegeneracy} and \ref{est penalty
term}, we obtain the following Theorem, with estimates that are
uniform in $p$.
\begin{theorem}\label{Estimates uniform in p}
    There exists a constant $K>0$, depending on dimension, $f$ and
    $\alpha$, but independent of $p$ such that for any solution
    $u_p$ of $\eqref{problem.p}$, there holds
    \begin{equation}\label{Uniform Lip}
        \| \nabla u_p \|_{L^\infty(\Omega)} \le K.
    \end{equation}
    Moreover, $u_p$ growth linearly uniform-in-$p$ away from the free boundary, that is,
    for a constant $\gamma > 0$ independent of $p$,
    \begin{equation}\label{Uniform Linear Growth}
        u_p(x) \ge \gamma \, \dist (x, \partial \{u_p > 0 \} ), \quad
        \forall x \in \{u_p > 0 \}.
    \end{equation}
    In addition, $u_p$ is uniformly strong nondegenerate, that is,
    for any fixed free boundary point $x_0 \in \partial \{ u_p > 0
    \}$,
    \begin{equation}\label{Uniform Nondeg}
        \sup\limits_{B_r(x_0)} u_p \ge \gamma r,
    \end{equation}
    where $\gamma > 0$ is independent of $p$.
\end{theorem}

\medskip

\noindent \textbf{Proof of Theorem \ref{Lip and Nondeg of Limit}.}
Notice that $\lim\limits_{p\to \infty} p^{1/p} = 1$. Passing the
limit as $p$ goes to infinity in \eqref{Uniform Lip},
\eqref{Uniform Linear Growth} and \eqref{Uniform Nondeg}, we prove
Theorem \ref{Lip and Nondeg of Limit}. \hfill $\square$

\medskip

\par

Theorem \ref{Estimates uniform in p} actually gives more
qualitative information than Theorem \ref{Lip and Nondeg of Limit}
itself. Indeed, with Theorem \ref{Estimates uniform in p} we can
address free boundary convergence issues. In what follows we prove
convergence of the free boundaries in the Hausdorff metric,
Theorem \ref{Hausdorff Convergence.intro}.

\medskip

\noindent \textbf{Proof of Theorem \ref{Hausdorff
Convergence.intro}.} For any set $A \subset \mathbb{R}^n$, and
$\varepsilon > 0$ fixed, let $\Gamma_\varepsilon(A)$ denote the
$\varepsilon$-neighborhood of $A$, that is,
$$
    \Gamma_\varepsilon(A) := \left \{ x \in \mathbb{R}^n \suchthat \dist
    (x, A) < \varepsilon \right \}.
$$
We have to show that given $\varepsilon > 0$, for $p \gg 1$,
depending on $\varepsilon > 0$, there hold
$$
    \begin{array}{c}
        \partial \{ u_p > 0 \} \subset \Gamma_\varepsilon \left ( \partial \{ u_\infty > 0
        \} \right ) \\
        \text{and} \\
        \partial \{ u_\infty > 0 \} \subset \Gamma_\varepsilon \left ( \partial \{ u_p > 0
        \} \right ).
    \end{array}
$$
Let $\xi$ be an arbitrary point on $\partial \{ u_p > 0 \}$ and
let us assume, for sake of contradiction, that $\xi \not \in
\Gamma_\varepsilon (\partial \{ u_\infty > 0 \}) $, that is,
$$
    \dist (\xi,\partial \{ u_\infty > 0 \}) \ge \varepsilon.
$$
If $u_\infty (\xi ) > 0$, then by linear growth, we would have
$$
    u_\infty(\xi) \ge \gamma \dist (\xi,\partial \{ u_\infty > 0
    \}) \ge \gamma \varepsilon,
$$
Thus, from uniform convergence, if $p \gg 1$, $u_p(\xi) \ge
\frac{2}{3} \gamma \varepsilon$, driving us to a contradiction. If
we assume $u_\infty(\xi) = 0$, then $u_\infty \big
|_{B_\varepsilon(\xi)} \equiv 0$. However, by strong
nondegeneracy, we know that
$$
    \sup\limits_{B_{\frac{\varepsilon}{2}}} u_p \ge \gamma_0
    \frac{\varepsilon}{2},
$$
and again it would drive us to a contraction on the uniform
convergence of $u_p$ to $u_\infty$. We have proven
$$
    \partial \{ u_p > 0 \} \subset \Gamma_\varepsilon \left ( \partial \{ u_\infty > 0 \}
     \right ).
$$
The other inclusion is proven similarly. \hfill $\square$

\medskip

\noindent \textbf{Proof of Theorem \ref{FBCond}.} Initially let us
recall some further facts from the $p$-Dirichlet minimization
problem \eqref{problem.p}. Recall that the free boundary $\partial
\{u_p > 0 \}$ is a $C^{1,\alpha}$ smooth surface up to a
$\mathcal{H}^{n-1}$ closed and negligible set (see
\cite{FMW06}, \cite{OT}, \cite{DP}). From the free boundary
condition $|\nabla u_p| = \lambda_p$, we deduce that
\begin{equation}\label{FBC Eq1}
    \lim\limits_{\begin{array}{l} x \to \partial \Omega_p \\
    x \in \Omega_p \end{array}} \dfrac{u_p(x)}{\dist(x,\partial
    \Omega_p)} = \lambda_p.
\end{equation}
Hereafter, $\Omega_p$ denotes the set of positivity of $u_p$. From
uniform convergence, $u_p \rightrightarrows u_\infty$, given a
point $x \in \Omega_\infty$, we may assume $x \in \Omega_p$ for
$p$ sufficiently large. Now, from the free boundary convergence
result, Theorem \ref{Hausdorff Convergence.intro}, there holds
\begin{equation}\label{FBC Eq2}
    \dist (x, \partial \Omega_\infty) = \dist (x, \partial
    \Omega_p) + \text{o}(1), \quad \text{ as } p \nearrow \infty.
\end{equation}
Here, $o(1)$ is an error that goes to zero as $p$ goes to
infinity. Thus, using once more the Hausdorff metric convergence
of the free boundary and uniform convergence of $u_p$ to
$u_\infty$, together with \eqref{FBC Eq1} and \eqref{FBC Eq2}, we
reach the following chain
$$
    \begin{array}{lll}
        \dfrac{u_\infty(x)}{\dist(x,\partial \Omega_\infty)} & = &
        \dfrac{u_p(x)}{\dist(x,\partial \Omega_p)} + \text{o}(1) \\
        &=& \lambda_p \dist(x,\partial \Omega_p) + \text{o}(\dist(x,\partial \Omega_p))
        + \text{o}(1)\\
        &=& \lambda_\infty \dist(x,\partial \Omega_\infty) +
        \text{o}(\dist(x,\partial \Omega_\infty)) + \text{o}(1)
    \end{array}
$$
Letting $p \to \infty$ the proof of Theorem \ref{FBCond} is
complete. \hfill $\square$

\bigskip
\bigskip

\noindent {\bf Acknowledgements.} J. D. Rossi partially supported
by MTM2004-02223, MEC, Spain, by UBA X066 and by CONICET,
Argentina. E. Teixeira's research has been partially supported by
CNPq-Brazil.

\bibliographystyle{amsplain, amsalpha}

\end{document}